\documentclass[12pt,a4paper]{article}
\usepackage[margin=1in]{geometry}
\usepackage{amssymb}
\usepackage{latexsym}
\usepackage{graphicx}
\usepackage{amsfonts}
\usepackage{amsthm}
\usepackage{enumerate}
\usepackage{amsmath}
\usepackage{lmodern}
\usepackage{cite}
\usepackage{MnSymbol}
\usepackage{amsmath,amscd}
\usepackage{pictexwd,dcpic}
\usepackage{parskip}

\theoremstyle{remark}

\theoremstyle{definition}

\author{Ewan Johnstone}
\title{Birationally rigid singular double quadrics and double cubics}

\begin{document}
\date{}
\maketitle

\begin{abstract}
In this paper it is shown that Fano double quadrics of index 1 and dimension at least 6 are birationally superrigid if the branch divisor has at most quadratic singularities of rank at least 6. Fano double cubics of index 1 and dimension at least 8 are birationally superrigid if the branch divisor has at most quadratic singularities of rank at least 8 and another minor condition of general position is satisfied. Hence, in the parameter spaces of these varieties the complement to the set of factorial and birationally superrigid varieties is of codimension at least $\binom{M-4}{2} + 1 $ and $\binom{M-6}{2} + 1 $ respectively.\\ Bibliography: 16 titles.
\end{abstract}

\section{Introduction}
Recall that a Fano double hypersurface of index 1 is defined as a projective algebraic variety $V$ equipped with a morphism $\sigma: V \to Q \subset \mathbb{P}^{M+1} $, where $M \geq 4 $, $\mathbb{P}^{M+1}$ is the complex projective space, $Q$ is an irreducible hypersurface of degree $m$ where $ 2 \leq m \leq M-2 $ and $\sigma$ is a double cover ramified over a divisor $W \subset Q$ which is cut out on $Q$ by a hypersurface $W^* \subset \mathbb{P}^{M+1}$ of degree $2M - 2m+2 $.
The variety $V$ can be realized as a complete intersection of codimension $2$ in the weighted projective space $ \mathbb{P}(1^{M+2},M-m+1) $, given by the equations

$$ f(x_0 , \dotsc, x_{M+1} ) = 0, \quad y^2 = g(x_0, \dotsc, x_{M+1}), $$

where $x_0, \dotsc, x_{M+1} $ have weight $1$ and $y$ has weight $M-m+1$. When $m=2$, we call $V$ a double quadric, and $m=3$ a double cubic. The hypersurface $Q$ is assumed to be non-singular throughout the paper. For shorthand we write $\mathbb{P}$ for $\mathbb{P}^{M+1}$.

Let ${\cal F}_m \subset \mathbb{P} (H^0(\mathbb{P},{\cal O}_{\mathbb{P}}(m))) $ be the open set of non-singular hypersurfaces $Q$ and let ${\cal G} = \mathbb{P} (H^0(\mathbb{P},{\cal O}_{\mathbb{P}}(2M -2m +2))) $ be the parameter space of hypersurfaces $W^*$. Then ${\cal I}_m = {\cal F}_m \times {\cal G } $ is a natural parameter space for double hypersurfaces of degree $2m$ with respect to the ample class $L = \sigma^*{\cal O}_Q (1) $. Let $S_m \subset {\cal I}_m $ be the set of pairs $(Q,W^*)$ such that the corresponding double hypersurface $V$ is a factorial variety with at most terminal singularities (and therefore, a Fano variety with $\mathop{\rm Pic} V = \mathbb{Z}K_V, K_V = -L$ ) which is birationally superrigid (we recall the definition of superrigidity in Section 2). The main result of this paper is the following claim.

{\bf Theorem 1.1. } {\it The complement  $\overline{{\cal I}_k \setminus S_k }$ has codimension at least $\binom{M-4}{2}+1 $ for $M \geq 6$ when $k =2$ and at least $\binom{M-6}{2}+1 $ for $M \geq 8 $ when $k =3 $.    }

This result is derived from two more explicit claims.

{\bf Theorem 1.2. } {\it Let $\sigma: V \to Q \subset \mathbb{P}$ be a double quadric ramified over $W = W^* \cap Q $. Assume that $M \geq 6$ and $W$ has at most quadratic singularities of rank at least $6$. Then $V$ is factorial and birationally superrigid.}
 
For the double cubic case, a small regularity condition is needed which we now introduce. Let $\sigma:V \to Q \subset \mathbb{P} $ be a double cubic, branched over $W = W^* \cap Q$. We say that $V$ satisfies the condition (*) if for any non-singular point $p \in W$ of the branch divisor in a system $(z_1, \dotsc, z_{M+1}) $ of affine coordinates with the origin at $p$ the hypersurfaces $Q$ is given by the equation

$$ 0 = q_1(z_*) + q_2(z_*) + \dotsc $$

where $q_i(z_*) $ are homogeneous of degree $i$ and $\left. q_2 \right \mid_{T_pW} \nequiv 0  $.

{\bf Theorem 1.3.} {\it Assume that $M \geq 8$, $W = W^* \cap Q$ has at most quadratic singularities of rank at least $8$ and $V$ satisfies the condition (*). Then $V$ is factorial and birationally superrigid. }

In the statements of Theorems 1.2 and 1.3 , the condition for the singularities to be quadratic of rank at least $r$ is understood in the sense of \cite{MR3229348}: for any point $o \in W$ and a system $z_1, \dotsc , z_M$ of analytic coordinates on $V$ with the origin at $o$, the branch divisor is given by an analytic equation

$$ 0 = z_1^2 + \dotsc + z_k^2 + h(z)$$ 
where $k \geq r$ and $h(z)$ stands for the terms of degree $3$ and higher, see \cite{MR3229348} for details.

The present paper continues a number of previous works. Superrigidity of generic (in particular, non-singular) double hypersurfaces was first shown in \cite{MR1777571}. Certain singular cases were investigated in \cite{MR2206637}, see also \cite{MR2195677}. Cyclic covers of degree $3$ and higher were studied in \cite{MR2508899}, triple spaces with isolated quadratic points in \cite{MR2108528}.

Double spaces of index $1$ with higher-dimensional singular locus were shown to be birationally superrigid in \cite{MR2761606}. Here we work in the style of \cite{MR3229348}, not only showing birational superrigidity of a certain class of Fano varieties but also estimating the codimension of the complement to the set of factorial and superrigid varieties. Such estimates are important due to applications to the theory of birational rigidity of Fano fibre spaces, see \cite{MR3397423}.

The structure of the paper is as follows; in section $2$ we recall the definition of birational superrigidity and prove Theorem 1.2. After that, we derive Theorem $1.1$ from Theorems $1.2$ and $1.3$ .
In section $3$ we prove Theorem $1.3$.

The author would like to thank Professor Aleksandr Pukhlikov for bringing this problem to their attention and for many useful conversations, comments and corrections.

\section{Birationally superrigid double quadrics}

In subsection 2.1 we recall the definitions of birational superrigidity and maximal singularity. In subsection 2.2 we prove Theorem 1.2. In subsection 2.3 we deduce Theorem 1.1 from Theorems 1.2 and 1.3.

{\bf 2.1. Birational superrigidity.}

Let $X$ be a factorial Fano variety with terminal singularities. For any effective divisor $D$ we define the {\it threshold of canonical adjunction} as

$$ c(D,X) = \mathop{\rm sup} \left \{ \frac{b}{a} \mid a \in \mathbb{Z}_+ \setminus 0 , b \in \mathbb{Z}_+, \lvert aD + bK_X \rvert \neq 0 \right \} .$$

{\bf Example 2.1.} 
If $V$ is a double quadric satisfying the assumptions of Theorem $1.2$ or a double cubic satisfying the assumptions of Theorem $1.3$ then it is easy to see that $V$ has at most quadratic (in particular, hypersurface) singularities and $\mathop{\rm codim}( \mathop{\rm Sing} V \subset V) \geq 6$ in the quadric case and $\mathop{\rm codim}( \mathop{\rm Sing} V \subset V) \geq 8$ in the cubic case. Therefore, as a corollary of Grothendieck's theorem on parafactoriality of local rings  (see \cite{MR1266175}), $V$ is factorial. By the Lefschetz hyperplane theorem we have $\mathop{\rm Pic} V = \mathbb{Z}L $ and by the adjunction formula, $K_V = -L$. So any effective divisor $D \sim nL , n \geq 1$ and then $c(D,V) = n$.
If $\Sigma$ is a linear system and $D \in \Sigma$, we set $c(\Sigma,X) = c(D,X) $.
The same definitions of the threshold of canonical adjunction (for a divisor and a linear system) apply to non-singular projective rationally connected varieties as in Chapter 2, p.38 of \cite{MR3060242}.

{\bf Definition 2.1.} 
For any mobile linear system $\Sigma$ on a factorial Fano variety $X$ with terminal singularities the {\it virtual threshold of canonical adjunction} is

$$ c_{\rm virt} (\Sigma) = \mathop{\rm inf}_{X^* \to X} \{ c(\Sigma^*, X^*) \} $$

where the infimum is taken over all birational morphisms $X^* \to X$ with $X^*$ smooth and projective. $\Sigma^*$ is the strict transform of $\Sigma$ on $X^*$. See Chapter 2, p.39 of \cite{MR3060242} for details.

{\bf Definition 2.2.}
A factorial Fano variety $X$ with at most terminal singularities is {\it birationally superrigid} if for any mobile linear system $\Sigma $ on $X$ the equality $$ c_{\rm virt}(\Sigma) = c(\Sigma,X) $$ holds.

It is well known that birationally superrigid varieties have a number of interesting geometric properties: they cannot be fibred into rationally connected varieties in a non-trivial way, their groups of birational self-maps $\mathop{\rm Bir} X$ and biregular automorphisms $\mathop{\rm Aut} X$ coincide, they are non-rational ect. See Chapter 2 of \cite{MR3060242} for more details.
Again, let $X$ be a factorial Fano variety with at most terminal singularities and $\mathop{\rm Pic} X = \mathbb{Z}K_X $.

{\bf Proposition 2.1. }
{\it Assume that $X$ is not birationally superrigid. Then there exists a mobile linear system $\Sigma \subset \lvert -nK_X \rvert $, a birational morphism $\phi: X^* \to X$ with $X^*$ smooth and projective and a $\phi$-exceptional divisor $E \subset X^*$ such that the Noether-Fano inequality  

$$ \mathop{\rm ord}\nolimits_E \phi^* \Sigma > na(E,X) $$

holds. }

\begin{proof}
See Chapter 2, Proposition 1.2, p.42 of\cite{MR3060242}. Although the case considered  there is non-singular the arguments work in the factorial terminal case also.
\end{proof}

The exceptional divisor $E$ (more precisely, the geometric discrete valuation $\mathop{\rm ord}\nolimits_E $) which satisfies the Noether-Fano inequality as stated in Proposition 2.1 is called a {\it maximal singularity} of the system $\Sigma$.

Thus, in order to prove that a given Fano variety is birationally superrigid one must show that there are no mobile linear systems with a maximal singularity on this variety.
An important particular case of a maximal singularity is a maximal subvariety. For $X, \Sigma, n \geq 1$ as above, we say that an irreducible subvariety $Y \subset X, Y \nsubset \mathop{\rm Sing} X $ of codimension at least $2$ is a {\it maximal subvariety} of a linear system $\Sigma$ if $ \mathop{\rm mult}\nolimits_Y \Sigma > n(\mathop{\rm codim}Y - 1) $.

{\bf 2.2. Double quadrics.}
Let us prove Theorem 1.2. Assume that $\Sigma \subset \lvert nL \rvert $ is a mobile linear system with a maximal singularity $E \subset V^* $, where $\phi: V^* \to V$ is a birational morphism from a non-singular projective variety $V^*$. Let $B = \phi(E) $ be the centre of $E$ on $V$. Assume first that $\mathop{\rm codim} B \geq 3$.

{\bf Lemma 2.1.}
{\it Let $D_1 , D_2 $ be generic divisors in $\Sigma$ and $Z = (D_1 \circ D_2 ) $ the cycle of their scheme-theoretic intersection. Then the inequality 

$$ \mathop{\rm mult}\nolimits_B Z > 4n^2 $$

holds.}

\begin{proof}
See Proposition $1$ in \cite{MR3229348} and Proposition 2.4 in \cite{2015arXiv}.

\end{proof}

Now as the linear system $\lvert L \rvert $ defines precisely the double cover $\sigma$, for any point $o \in B$ we get the inequality

$$ \mathop{\rm mult}\nolimits_o Z \leq \mathop{\rm deg}\nolimits_L Z = 4n^2 $$

We conclude that the case when $\mathop{\rm codim} B \geq 3$ is impossible.

We now consider the case when $B$ is a subvariety of codimension 2.
Then it is a maximal subvariety of the linear system $\Sigma$ and the inequality $\mathop{\rm mult}\nolimits_B \Sigma > n$ holds. We define $V_H = V \cap H$, where $H$ is a generic linear subvariety of dimension $6$. Since the codimension of the singular set of $V$ is at least $6$ we conclude that $V_H$ is non-singular. We define $B_H = B \cap H$ and note that it satisfies the same inequality with respect to $\Sigma_H = \left. \Sigma \right \mid_{V_H}$.   Since $\mathop{\rm dim}V_H \geq 5$ by the Lefschetz hyperplane theorem  $B_H$ is numerically equivalent to a multiple $kL^2$. Set $ \nu = \mathop{\rm mult}\nolimits_{B_H} \Sigma_H > n $. Then for the cycle $Z$ defined in Lemma 2.1 we get $Z \sim n^2L^2 $ and 
$$Z = \nu^2B_H + Z_1 $$

with $Z_1$ an effective cycle. Comparing the $L$-degrees on the left and right hand sides, we get a contradiction. Thus concluding the proof of Theorem 1.2.

{\bf Remark 2.1.}
More subtle arguments give a proof of Theorem $1.2$ for $M=4$ and $5$ under the slightly weaker assumption that the quadratic singularities are of rank at least $4$. Only the arguments for the case $\mathop{\rm codim} B = 2$ need to be modified in a way similar to \cite{MR936532}. 

{\bf 2.3. Codimension in the parameter space.}
Let us prove Theorem 1.1, assuming Theorems 1.2 and 1.3. First, we consider double quadrics.
Fix a non-singular hypersurface $Q \in {\cal F}_m$. Given the claim of Theorem 1.2, it is sufficient to show that the set of hypersurfaces $W^* \in {\cal G} $ violating the assumptions of Theorem 1.2 is of codimension at least $\binom{M-4}{2} +1 $ in ${\cal G}$. This is easily done by a standard dimension count.
Fix a point $p \in Q$. The hypersurface $Q$ is given by an equation

$$ 0 = q_1 + q_2 + \dotsc $$

in some system of affine coordinates with origin at $p$. The hypersurface $W^*$ is given by an equation

$$ 0 = w_1 + w_2 + \dotsc $$

where $w_i$ are homogeneous of degree $i$ (we assume that $p \in W$ - otherwise the case is trivial). Violation of the assumptions of Theorem 1.2 at $p$ means that $w_1 = \lambda q_1$ for some $\lambda \in \mathbb{C}$ and  $\left. w_2 \right \mid_{q_1=0} $ is a quadratic form of rank at most $5$. This imposes

$$ \frac{(M-5)(M-4)}{2} + M $$

independent conditions on the coefficients of the polynomial $g$. Now the standard dimension count (the point $p$ varies in the intersection $W = W^* \cap Q$) completes the proof of Theorem 1.1 in the case of double quadrics.
In the case of double cubics we have almost identically the same estimates for the condition on the singularities of $W$.
However, we also have the condition (*). Again, we consider first the violation of (*) at a fixed point $p \in W$. The condition

 $$\left. q_2 \right \mid_{T_p W} \equiv 0 $$

for fixed $q_1, w_1$ (in the notations above) imposes $\frac{M(M-1)}{2}$ independent conditions on $q_2$. Therefore, the set of pairs $(f,g)$, such that in a system of affine coordinates with the origin at some point $p \in \{ f=g=0 \} $ the condition (*) is violated is of codimension at least $\frac{M(M-1)}{2} - M = \frac{M(M-3)}{2} $ in the parameter space. As this number is higher than $\binom{M-6}{2} + 1 $, the proof of Theorem 1.1 is complete.

\section{Birationally superrigid double cubics}

In this section we prove Theorem 1.3. First we exclude maximal singularities the centre of which is not contained in the singular locus $\mathop{\rm Sing}V$ (Subsection 3.1). We then exclude maximal singularities with the centre inside $\mathop{\rm Sing} V$ (Subsections 3.2 and 3.3). 

{\bf 3.1. Maximal singularities outside the singular locus.} 
Recall that we work with a double cover $\sigma: V \to Q$, where $Q \subset \mathbb{P}$ is a non-singular cubic hypersurface, $\sigma$ is branched over $W = W^* \cap Q$ and the assumptions of Theorem 1.3 holds. In particular, $\mathop{\rm Pic} V = \mathbb{Z}L $ where $L$ is the $\sigma$-pullback of the hyper-plane section of $Q$. Assume that $V$ is not birationally superrigid. Then there is a mobile linear system $\Sigma \subset \lvert nL \rvert $ with a maximal singularity $E \subset V^*$, where $\phi:V^* \to V$ is a birational morphism from a non-singular projective $V^*$.
Assume first that $B = \phi(E) \not\subset \mathop{\rm Sing}V$. Note that this implies $ \sigma(B) \not\subset \mathop{\rm Sing}W $.

If $\mathop{\rm codim} B = 2$, we argue as for double quadrics (Subsection 2.2) and come to a contradiction. if $\mathop{\rm codim}B \geq 3$, then the self-intersection $Z = (D_1 \circ D_2)$ of the system $\Sigma$ (where $D_1, D_2$ are general divisors) satisfy the $4n^2$-inequality $\mathop{\rm mult}\nolimits_B Z > 4n^2 $. If $\sigma(B) \subset W$, then we argue as in the end of section 3 of \cite{MR1777571}.
Take a point $p \in \sigma(B) \setminus \mathop{\rm Sing} W$. By the condition (*), the irreducible subvariety $\sigma^{-1}(T_p W \cap Q) \sim L^2 $ has multiplicity precisely $4$ at $o = \sigma^{-1}(p)$. Therefore, there exists an irreducible component $Y$ of the cycle $Z$, such that $Y \sim lL^2 , \mathop{\rm mult}\nolimits_o Y  > 4l $  and $Y \neq \sigma^{-1} (T_p W \cap Q)$. So $Y$ is not contained in both divisors $\sigma^{-1} (T_p Q \cap Q)$ and $\sigma^{-1} (T_p W^* \cap Q)$.

Intersecting $Y$ with one which does not contain $Y$, we obtain an effective cycle $Y^* \sim lL^3 $ of codimension $3$ such that $\mathop{\rm mult}\nolimits_o Y^* > 8l$. As $\mathop{\rm deg}\nolimits_L Y^* = 6l$, we obtain a contradiction.

The same argument works in the case of $B$ such that $\sigma(B) \not\subset  W$ but also $\sigma(B) \cap W \not\subset \mathop{\rm Sing} W$, provided we take $p$ such that $p \in \sigma(B) \cap W$, else $\sigma^{-1}(p)$ is not well defined. Therefore, the last case to consider is a maximal singularity such that its centre $B$ satisfies the property $\sigma(B) \not\subset W$ and  $\sigma(B) \cap W \subset \mathop{\rm Sing} W$.

In that case $\mathop{\rm codim} B \geq 6 $. We take a general point $o \in B$ such that $\bar{o} = \sigma(o) \not\in W$. Let $\psi: V^+ \to V$ be the blow up of the point $o$ and let $E = \psi^{-1}(o)$ be the exceptional divisor. By the $8n^2$-inequality (See Chapter 2, Theorem 4.1, p74 of \cite{MR3060242}, \cite{MR2533620}), there is a linear subspace $\Psi \subset E \cong \mathbb{P}^{M-1} $ of codimension $2$, such that the self-intersection $Z$ of $\Sigma$ satisfies the property

$$ \mathop{\rm mult}\nolimits_o Z + \mathop{\rm mult}\nolimits_{\Psi} Z^+ > 8n^2 $$

$Z^+$ being its strict transform on $V^+$. As $\sigma$ gives an isomorphism of local rings ${\cal O}_{o,V} $ and ${\cal O}_{\bar{o},Q} $, we can find a hyperplane $H \subset \mathbb{P} $, such that $\sigma^{-1} (H_Q)^+ $ contains $\Psi$, where $H_Q = H \cap Q$ is the corresponding hyperplane section. As such hypersurfaces form a 2-dimensional linear system, we may assume that $\sigma^{-1} (H_Q) $ contains none of the components of $Z$, so the cycle $( Z \circ H_Q)$ is well defined. Obviously, $(Z \circ H_Q) \sim n^2L^3 $, so $\mathop{\rm deg}(Z \circ H_Q) = 6n^2 $. On the other hand, 
$$ \mathop{\rm mult}\nolimits_o (Z \circ H_Q) \geq \mathop{\rm mult}\nolimits_o Z + \mathop{\rm mult}\nolimits_{\Psi}Z^+ > 8n^2  $$ 
a contradiction. \\
We have inspected all options for the case $B \not\subset \mathop{\rm Sing} V$. Therefore we may assume that $B \subset \mathop{\rm Sing} V$ and so $\sigma(B) \subset \mathop{\rm Sing}W$.

{\bf 3.2. Inversion of adjunction.}
Fix a general point $o \in B$. As the singularities of $W$  are quadratic of rank at least $8$, the singularities of $V$ are quadratic of rank at least $9$. Near the point $o$ we may consider the germ $o \in V$ analytically as a germ of a hypersurface in $\mathbb{C}^{M+1} $. Let $ X \ni o $ be the section of the germ $o \in V$ by a generic $5$-plane containing the point $o$. Then $o \in X$ is a germ of an isolated quadratic singularity of rank $5$.
Let $\pi: V^+ \to V$ and $\pi_X: X^+ \to X$ be the blow ups of the point $o$, and $E = \pi^{-1}(o), E_X = \pi^{-1}_X(o) $ the exceptional divisors. In an obvious sense, the non-singular 3-dimensional quadric $E_X$ is the section of the quadric $E$ by a generic 4-plane. For a general divisor $D \in \Sigma$ set

$$ \pi^*D = D^+ + \nu E $$

By $D_X$ we denote the restriction of $D$ onto $X$ and by $D^+_X$ its strict transform on $X^+$. By inversion of adjunction, the pair $(X, \frac{1}{n}D_X)$ is not log canonical at $o$.
By our assumption about $B$ and what was shown in Subsection 3.1, the point $o$ is an isolated centre of a non-canonical singularity of the pair $(X,\frac{1}{n}D_X) $.

{\bf Proposition 3.1.}
{\it The multiplicity $\nu$ satisfies the inequalities

$$ n < \nu \leq \sqrt{3}n $$
}

\begin{proof}

For the inequality $v > n$, see Chapter 7, Proposition 2.7, p.308 of \cite{MR3060242}. For the inequality $\nu \leq \sqrt{3}n$, we note that $\mathop{\rm deg}\nolimits_L Z = 6n^2 \geq \mathop{\rm mult}\nolimits_o Z \geq 2\nu^2 $.
\end{proof}

In particular, $\nu < 2n$, which implies that the pair $(X^+, \frac{1}{n} D_X^+ ) $ is non log-canonical and the centers of log-canonical singularities are contained in $E_X$. Moreover, by the connectedness principal (Chapter 17 of \cite{MR1225842})  the union of centres of all non log-canonical singularities is a connected subset of $E_X$. We therefore have 3 cases to consider:

\begin{enumerate}
 \item $(X^+, \frac{1}{n}D_X^+)$ is non log-canonical at a surface in $E_X$. 
 \item $(X^+, \frac{1}{n}D^+_X) $ is non log-canonical at a curve in $E_X$. 
\item $(X^+, \frac{1}{n}D^+_X)$ is non log-canonical at a point $p_X \in E_X $.
\end{enumerate}

{\bf Lemma 3.1.} {\it Case 1 cannot occur.}
\begin{proof}
In this case the pair $(X^+, \frac{1}{n} \Sigma^+_X)$ is non log-canonical at an irreducible divisor $R_X \subset E_X$, which is a section of an irreducible divisor $R \subset E$. Therefore, the ${\cal O}_E(1) $-degree of $R$ is at least $2$.
We now compute as in the proof of Theorem 4.1 in Chapter 2 of \cite{MR3060242} (see also \cite{MR2195677}, \cite{MR2533620}). Using the fact that $(X^+, \frac{1}{n} \Sigma^+_X)$ is non log-canonical,
for generic divisors $D_1 , D_2 \in \Sigma $ we obtain the inequality

$$ \mathop{\rm mult}\nolimits_R (D^+_1 \circ D^+_2 ) > 4n^2 $$

so that $\mathop{\rm mult}_o Z \geq 2\nu^2 + 8n^2 > 10n^2$, a contradiction.
\end{proof}

{\bf Lemma 3.2.} {\it Case 3 cannot occur.}

\begin{proof}
We note that $p_X = S \cap E_X $ where $S$ must be a linear subspace of codimension $3$ on $E$, as $E_X$ is a generic section of $E$. Since a quadric of rank at least $9$ does not contain a linear subspace of codimension $3$, we arrive at a contradiction.
\end{proof}
Therefore the only case left is case 2: $(X^+, \frac{1}{n}D^+_X) $ is non log-canonical at an irreducible curve $Y_X$ which is a section of an irreducible subvariety $Y \subset E$ of codimension $2$. 

{\bf 3.3. Centre at an irreducible curve.}
Write the self-intersection of the linear system $\Sigma_X^+$ as

$$ (\left. D_1^+ \right \mid _X \circ \left. D_2^+ \right \mid _X ) = Z^+_X + Z^*_X $$

where $Z^*_X$ is an effective divisor on $E_X$, which is the restriction onto $E_X$ of an effective divisor $Z^*$ on $E$. We have the $4n^2$-inequality

$$ \mathop{\rm mult}\nolimits_{Y_X} Z^+_X + \mathop{\rm mult}\nolimits_{Y_X} Z^*_X > 4n^2 $$

so that $\mathop{\rm mult}\nolimits_Y Z^+ + \mathop{\rm mult}\nolimits_{Y} Z^* > 4n^2 $.
We first give a bound on the multiplicity of $\mathop{\rm mult}\nolimits_Y Z^* $.

{\bf Lemma 3.3.} {\it The inequality $\mathop{\rm mult}\nolimits_Y Z^* = \mathop{\rm mult}\nolimits_{Y_X} Z^*_X \leq 2n^2 $ holds.}

\begin{proof}
Assume the converse: $\mathop{\rm mult}\nolimits_Y Z^* > 2n^2 $. As $Z^*_X$ is an effective divisor on a non-singular quadric, by the cone technique (See Chapter 2, Proposition 2.3 of \cite{MR3060242} ) we have

$$ Z^* \sim \alpha H_E $$

where $ \alpha > 2n^2 $, $H_E$ is the hyperplane section of $E$. However, $ \mathop{\rm mult}\nolimits_o Z > 2\nu^2 + 4n^2 > 6n^2 $, a contradiction.
\end{proof}

From this we conclude that $\mathop{\rm mult}\nolimits_Y Z^+ = \mathop{\rm mult}\nolimits_{Y_X} Z^+_X > 2n^2 $.
Our next observation is that $Y$ must be a section of $E$ by a linear subspace of codimension 2 (recall that $E$ is a quadric with a natural embedding in projective space). Indeed, since $E$ is a quadric of rank at least $9$, $Y$ is numerically equivalent to $\beta H^2_E$ for some $\beta \geq 1$.

{\bf Lemma 3.4.} {\it The equality $\beta=1 $ holds.}
\begin{proof}
Assume the converse: $\beta \geq 2$. Then $\mathop{\rm deg} Y \geq 4$ so $\mathop{\rm mult}\nolimits_o Z = \mathop{\rm deg}(Z^+ \circ E) \geq \mathop{\rm mult}\nolimits_Y Z^+ \mathop{\rm deg} Y   > 8n^2 $. A contradiction.
\end{proof}

Therefore $\beta = 1$ and $Y$ is a section of $E$ by a linear subspace as claimed.

Set $\bar{o} = \sigma(o) \in Q$. Let $\pi_Q: Q^+ \to Q$ be the blow up of $\bar{o}$,  $ \overline{E} = \pi_Q^{-1}(\bar{o}) $ the exceptional divisor, $\overline{E} \cong \mathbb{P}^{M-1}$. Note that $\sigma: V \to Q$ extends to $\sigma^+: V^+ \to Q^+ $, where $\left. \sigma^+ \right \mid_E : E \to \overline{E} $ is a double cover branched over the quadric $W^+ \cap \overline{E} $. We now have two cases: 
\begin{itemize}
\item $\sigma^+(Y) $ is a linear subspace in $\overline{E}$ of codimension 2 and $\left. \sigma^+ \right \mid _Y $ is a double cover
\item $\sigma^+(Y) = \overline{Y}$ is a quadric in a hyperplane in $\overline{E}$ and $\left. \sigma^+ \right \mid _Y$ is birational.
\end{itemize}

The first case can be excluded by the same arguments used in subsection 3.1, where the $8n^2$-inequality was applied. We now consider the second case. For a quadric hypersurface cone $\Lambda \subset \mathbb{P}$ with vertex $\bar{o}$, set
$ \Lambda_Q = \left. \Lambda \right \mid_Q  = \Lambda \cap Q$. Take a general cone $\Lambda$ such that $\Lambda^+_Q \cap \overline{E} $ contains $\overline{Y}$, so that $\sigma^{-1}(\Lambda_Q)^+ $ contains $Y$. By generality, $\sigma^{-1}(\Lambda_Q) $ contains none of the components of $Z$ so the cycle

$$ Z_{\Lambda} = (Z \circ \sigma^{-1}(\Lambda_Q)) $$

is well defined. Its $L$-degree is $12n^2$. Now by standard intersection theory (See Chapters 1-3, \cite{MR1644323} or Chapter 2, Lemma 2.2 of \cite{MR3060242} and \cite{MR1798981} for details) we obtain

$$ \mathop{\rm mult}\nolimits_o Z_{\Lambda} \geq 2\mathop{\rm mult}\nolimits_o Z + 2\mathop{\rm mult}\nolimits_Y Z^+ > 8n^2 + 4n^2 = 12n^2 $$

which gives the final contradiction. This concludes the proof of Theorem 1.3.

\bibliographystyle{alpha}
\bibliography{hyper}{}

\end{document}